\documentclass[thmsa,12pt]{article}
\usepackage{amsfonts}
\usepackage{amssymb}
\usepackage{amsmath}

\input{tcilatex}
\begin{document}

\begin{center}
{\Large Common Fixed Point Results for Family of Generalized Multivalued
F-contraction Mappings in Ordered Metric Spaces}{\large \medskip }

\textbf{Talat Nazir}$^{a,b}$\textbf{\ and Sergei Silvestrov}$^{b}$

{\small $^{(a)}$\textit{Division of Applied Mathematics, School of
Education, Culture and Communication, }}

{\small \textit{M\"{a}lardalen University, 72123 V\"{a}ster\aa s, Sweden.}}

$^{(b)}${\small \textit{Department of Mathematics, COMSATS Institute of
Information Technology,}}

{\small \textit{22060 Abbottabad, Pakistan.}}

E-mail{\small : talat@ciit.net.pk, sergei.silvestrov@mdh.se}\bigskip
\end{center}

\noindent
--------------------------------------------------------------------------------------------

\noindent \textit{Abstract} \_\_\_ In this paper, we study the existence of
common fixed points of family of multivalued mappings satisfying generalized 
$F$-contractive conditions in ordered metric spaces. These results establish
some of the general common fixed point theorems for family of multivalued
maps.

\noindent \textbf{---------------------------------------------}


\noindent \textit{2010 \ Mathematics Subject Classification: }\texttt{%
47H10,54H25, 54C60,46B40}.

\noindent \textbf{---------------------------------------------}

\noindent \textbf{1. \ \ INTRODUCTION AND PRELIMINARIES}

Markin \cite{Markin} initiated the study of fixed points for multivalued
nonexpansive and contractive maps. Later, a useful and interesting fixed
point theory for such maps was developed. Later, a rich and interesting
fixed point theory for such multivalued maps was developed; see, for
instance \cite{AlThag07, Ber07, Ciric09, Jung98, Kam, Lat, Lazar, Mot, Nad,
Rhoad, SP12}. The theory of multivalued maps has various applications in
convex optimization, dynamical systems, commutative algebra, differential
equations and economics. Recently, Wardowski \cite{Wardowski} introduced a
new contraction called $F-$contraction and proved a fixed point result as a
generalization of the Banach contraction principle. Abbas et al. \cite{Abb3}
obtained common fixed point results by employed the $F-$ contraction
condition. Further in this direction, Abbas et al. \cite{Abb2} introduced a
notion of generalized $F-$ contraction mapping and employed there results to
obtain a fixed point of a generalized nonexpansive mappings on star shaped
subsets of normed linear spaces. Minak et al. \cite{Minak} proved some fixed
point results for Ciric type generalized $F-$ contractions on complete
metric spaces. Recently, \cite{Acar} established some fixed point results
for multivalued $F-$contraction maps on complete metrics spaces.

The aim of this paper is to prove common fixed points theorems for family of
multivalued generalized $F$-contraction mappings without using any
commutativity condition in partially ordered metric space. These results
extend and unify various comparable results in the literature (\cite{Kannan}%
, \cite{LatifBeg}, \cite{Rus} and \cite{Sgr}).

We begin with some basic known definitions and results which will be used in
the sequel. Throughout this article, $%
\mathbb{N}
,$ $%
\mathbb{R}
^{+},$\ $%
\mathbb{R}
$\ denote the set of natural numbers, the set of positive real numbers and
the set of real numbers, respectively.

Let $\digamma $\ be the collection of all mappings $F:%
\mathbb{R}
^{+}\rightarrow 
\mathbb{R}
$ that satisfy the following conditions:

\begin{itemize}
\item[($F_{1}$)] $F$ is strictly increasing, that is, for all $a,b\in 
\mathbb{R}
^{+}$ \ such that $a<b$ implies that $F(a)<F(b)$.

\item[($F_{2}$)] For every sequence $\{a_{n}\}$ of positive real numbers, $%
\lim\limits_{n\rightarrow \infty }a_{n}=0$ and $\lim\limits_{n\rightarrow
\infty }F\left( a_{n}\right) =-\infty $ are equivalent.

\item[($F_{3}$)] There exists $\lambda \in \left( 0,1\right) $ such that $%
\lim\limits_{a\rightarrow 0^{+}}a^{\lambda }F(\lambda )=0$.
\end{itemize}

\noindent \textbf{Definition 1.1.} \cite{Wardowski} Let $(X,d)$ be a metric
space and $F\in \digamma $. A mapping $f:X\rightarrow X$ is said to be an $%
F- $contraction on $X$ if there exists $\tau >0$ such that%
\begin{equation*}
d(fx,fy)>0\text{ implies that }\tau +F\left( d(fx,fy)\right) \leq F\left(
d(x,y)\right) 
\end{equation*}%
for all $x,y\in X$.

Wardowski \cite{Wardowski} gave the following result.

\noindent \textbf{Theorem 1.2.} \ \ Let $(X,d)$ be a complete metric space
and mapping $f:X\rightarrow X$ be and $F-$contraction. Then there exists a
unique $x$ in $X$ such that $x=fx$. Moreover, for any $x_{0}\in X$, the
iterative sequence $x_{n+1}=f\left( x_{n}\right) $ converges to $x$%
.\smallskip

\noindent Kannan \cite{Kannan} has proved a fixed point theorem for a single
valued self mapping $T$ of a metric space $X$ satisfying the property%
\begin{equation*}
d(Tx,Ty)\leq h\{d(x,Tx)+d(y,Ty)\} 
\end{equation*}%
for all $x,y$ in $X$ and for a fixed $h\in \lbrack 0,\frac{1}{2}).$

Ciric \cite{Ciric} considered a mapping $T:X\rightarrow X$ satisfying the
following contractive condition:%
\begin{equation*}
d(Tx,Ty)\leq q\max \{d(x,y),d(x,Tx),d(y,Ty),d(x,Ty),d(y,Tx)\} 
\end{equation*}%
where $q\in \lbrack 0,1)$. He proved the existence of a fixed point when $X$
is a $T-$ orbitally complete metric space.\smallskip

Latif and Beg \cite{LatifBeg} extended Kannan mapping for multivalued
mapping and introduced the notion of a $K-$ multivalued mapping. Rus \cite%
{Rus} coined the term $R-$ multivalued mapping, which is a generalization of
a $K-$ multivalued mapping (see also, \cite{AAN12}). Abbas and Rhoades \cite%
{AbbasRhoades} studied common fixed point problems for multivalued mappings
and introduced the notion of generalized $R-$ multivalued mappings which in
turn\ generalizes $R-$ multivalued mappings.\smallskip

\noindent Let $(X,d)$ be a metric space. Denote by $P(X)$ be the family of
all nonempty subsets of $X$, and by $P_{cl}\left( X\right) $ the family of
all nonempty closed subsets of $X.$

\noindent A point $x$ in $X$ is called fixed point of a multivalued mapping $%
T:X\rightarrow P_{cl}(X)$ provided $x\in Tx.$ The collection of all fixed
point of $T$ is denoted by $Fix(T).$

Recall that, a map $T:X\rightarrow P_{cl}\left( X\right) $ is said to be
upper semicontinuous, if for $x_{n}\in X$ and $y_{n}\in Tx_{n}$ with $%
x_{n}\rightarrow x_{0}$ and $y_{n}\rightarrow y_{0}$, implies $y_{0}\in
Tx_{0}$.\smallskip

\noindent \textbf{Definition 1.3.} \ \ Let $X$ be a nonempty set. Then ($%
X,d,\preceq $) is called partially ordered metric space if and only if $d$
is a metric on a partially ordered set ($X,\preceq $).

We define $\Delta _{1},\Delta _{2}\subseteq X\times X$ as follows:%
\begin{eqnarray*}
\Delta _{1} &=&\{(x,y)\in X\times X\shortmid x\preceq y\}\text{ and} \\
\Delta _{2} &=&\{(x,y)\in X\times X\shortmid x\prec y\}.
\end{eqnarray*}

\noindent \textbf{Definition 1.4.} \ \ A subset $\Gamma $ of a partially
ordered set $X$ is said to be well-ordered if every two elements of $\Gamma $
are comparable.

\section{COMMON\ FIXED POINT RESULTS}

\noindent In this section, we obtain common fixed point theorems for family
of multivalued mappings. We begin with the following result.

\noindent \textbf{Theorem 2.1.} \ \ Let $(X,d,\preceq )$ be an ordered
complete metric space and $\{T_{i}\}_{i=1}^{m}:X\rightarrow P_{cl}(X)$ be
family of multivalued mappings. Suppose that for every $(x,y)\in \Delta _{1}$
and $u_{x}\in T_{i}(x),$ there exists $u_{y}\in T_{i+1}(y)$ for $i\in
\{1,2,...,m\}$ (with $T_{m+1}=T_{1}$ by convention) such that, $%
(u_{x},u_{y})\in \Delta _{2}$ implies%
\begin{equation}
\tau +F\left( d(u_{x},u_{y})\right) \leq F(M(x,y;u_{x},u_{y})),  \tag{2.1}
\end{equation}%
where $\tau $\ is a positive real number and%
\begin{equation*}
M(x,y;u_{x},u_{y})=\max \{d(x,y),d(x,u_{x}),d(y,u_{y}),\dfrac{d\left(
x,u_{y}\right) +d\left( y,u_{x}\right) }{2}\}. 
\end{equation*}%
Then the following statements hold:

\begin{enumerate}
\item[(1)] $Fix(T_{i})\neq \emptyset $ for any $i\in \{1,2,...,m\}$ if and
only if $Fix(T_{1})=Fix(T_{2})=...=Fix(T_{m})\neq \emptyset .$

\item[(2)] $Fix(T_{1})=Fix(T_{2})=...=Fix(T_{m})\neq \emptyset $ provided
that any one $T_{i}$ for $i\in \{1,2,...,m\}$\ is upper semicontinuous.

\item[(3)] $\cap _{i=1}^{m}Fix(T_{i})$ is well-ordered if and only if $\cap
_{i=1}^{m}Fix(T_{i})$ is singleton set.
\end{enumerate}

\noindent \textbf{Proof.} \ \ To prove (1), let $x^{\ast }\in T_{k}(x^{\ast
})$ for some $k\in \{1,2,...,m\}.$ Assume that $x^{\ast }\notin
T_{k+1}\left( x^{\ast }\right) ,$ then there exists an $x\in T_{k+1}\left(
x^{\ast }\right) $ with $\left( x^{\ast },x\right) \in \Delta _{2}\ $such
that%
\begin{equation*}
\tau +F\left( d(x^{\ast },x)\right) \leq F(M(x^{\ast },x^{\ast };x^{\ast
},x)), 
\end{equation*}%
where%
\begin{eqnarray*}
M(x^{\ast },x^{\ast };x^{\ast },x) &=&\max \{d(x^{\ast },x^{\ast
}),d(x^{\ast },x^{\ast }),d(x,x^{\ast }),\dfrac{d(x^{\ast },x)+d(x^{\ast
},x^{\ast })}{2}\} \\
&=&d(x,x^{\ast }),
\end{eqnarray*}%
implies that%
\begin{equation*}
\tau +F\left( d(x^{\ast },x)\right) \leq F(d(x^{\ast },x)), 
\end{equation*}%
a contradiction as $\tau >0$. Thus $x^{\ast }=x$. Thus $x^{\ast }\in
T_{k+1}\left( x^{\ast }\right) $ and so $Fix(T_{k})\subseteq Fix(T_{k+1}).$
Similarly, we obtain that\ $Fix(T_{k+1})\subseteq Fix(T_{k+2})$ and
continuing this way, we get $Fix(T_{1})=Fix(T_{2})=...Fix(T_{k}).$ The
converse is straightforward.\smallskip

\noindent To prove (2), suppose that $x_{0}$ is an arbitrary point of $X.$
If $x_{0}\in T_{k_{0}}\left( x_{0}\right) $ for any $k_{0}\in \{1,2,...,m\},$
then by using (1), the proof is finishes. So we assume that $x_{0}\notin
T_{k_{0}}\left( x_{0}\right) $ for any $k_{0}\in \{1,2,...,m\}.$ Now for $%
i\in \{1,2,...,m\}$, if $x_{1}\in T_{i}(x_{0}),$ then there exists $x_{2}\in
T_{i+1}(x_{1})$ with $(x_{1},x_{2})\in \Delta _{2}$ such that%
\begin{equation*}
\tau +F\left( d(x_{1},x_{2})\right) \leq F(M(x_{0},x_{1};x_{1},x_{2})), 
\end{equation*}%
where%
\begin{eqnarray*}
M(x_{0},x_{1};x_{1},x_{2}) &=&\max
\{d(x_{0},x_{1}),d(x_{0},x_{1}),d(x_{1},x_{2}),\dfrac{%
d(x_{0},x_{2})+d(x_{1},x_{1})}{2}\} \\
&=&\max \{d(x_{0},x_{1}),d(x_{1},x_{2}),\dfrac{d(x_{0},x_{2})}{2}\} \\
&=&\max \{d(x_{0},x_{1}),d(x_{1},x_{2})\}.
\end{eqnarray*}%
Now, if $M(x_{0},x_{1};x_{1},x_{2})=d(x_{1},x_{2})$ then%
\begin{equation*}
\tau +F\left( d(x_{1},x_{2})\right) \leq F(d(x_{1},x_{2})), 
\end{equation*}%
a contradiction as $\tau >0$. Therefore $%
M(x_{0},x_{1};x_{1},x_{2})=d(x_{0},x_{1})$ and we have%
\begin{equation*}
\tau +F\left( d(x_{1},x_{2})\right) \leq F\left( d(x_{0},x_{1})\right) . 
\end{equation*}%
Next for this $x_{2}\in T_{i+1}\left( x_{1}\right) ,$ there exists $x_{3}\in
T_{i+2}(x_{2})$ with $\left( x_{2},x_{3}\right) \in \Delta _{2}$ such that%
\begin{equation*}
\tau +F\left( d(x_{2},x_{3})\right) \leq F(M(x_{1},x_{2};x_{2},x_{3})), 
\end{equation*}%
where%
\begin{eqnarray*}
M(x_{1},x_{2};x_{2},x_{3}) &=&\max
\{d(x_{1},x_{2}),d(x_{1},x_{2}),d(x_{2},x_{3}),\dfrac{%
d(x_{1},x_{3})+d(x_{2},x_{2})}{2}\} \\
&=&\max \{d(x_{1},x_{2}),d(x_{2},x_{3})\}.
\end{eqnarray*}%
Now, if $M(x_{1},x_{2};x_{2},x_{3})=d(x_{2},x_{3})$ then%
\begin{equation*}
\tau +F\left( d(x_{2},x_{3})\right) \leq F(d(x_{2},x_{3})), 
\end{equation*}%
a contradiction as $\tau >0$. Therefore $%
M(x_{1},x_{2};x_{2},x_{3})=d(x_{1},x_{2})$ and we have%
\begin{equation*}
\tau +F\left( d(x_{2},x_{3})\right) \leq F\left( d(x_{1},x_{2})\right) . 
\end{equation*}%
Continuing this process, for $x_{2n}\in T_{i}(x_{2n-1})$, there exist $%
x_{2n+1}\in T_{i+1}\left( x_{2n}\right) $ with $\left(
x_{2n},x_{2n+1}\right) \in \Delta _{2}$ such that%
\begin{equation*}
\tau +F\left( d(x_{2n},x_{2n+1})\right) \leq F\left(
M(x_{2n-1},x_{2n};x_{2n},x_{2n+1})\right) , 
\end{equation*}%
where%
\begin{eqnarray*}
M(x_{2n-1},x_{2n};x_{2n},x_{2n+1}) &=&\max
\{d(x_{2n-1},x_{2n}),d(x_{2n-1},x_{2n}),d(x_{2n},x_{2n+1}), \\
&&\dfrac{d(x_{2n-1},x_{2n+1})+d(x_{2n},x_{2n})}{2}\} \\
&=&\{d(x_{2n-1},x_{2n}),d(x_{2n},x_{2n+1}),\dfrac{d(x_{2n-1},x_{2n+1})}{2}\}
\\
&\leq &d(x_{2n-1},x_{2n}),
\end{eqnarray*}%
that is,%
\begin{equation*}
\tau +F\left( d(x_{2n},x_{2n+1})\right) \leq F\left(
d(x_{2n-1},x_{2n})\right) . 
\end{equation*}%
Similarly, for $x_{2n+1}\in T_{i+1}(x_{2n})$, there exist $x_{2n+2}\in
T_{i+2}\left( x_{2n+1}\right) $ such that for $\left(
x_{2n+1},x_{2n+2}\right) \in \Delta _{2}$ implies%
\begin{equation*}
\tau +F\left( d(x_{2n+1},x_{2n+2})\right) \leq F\left(
d(x_{2n},x_{2n+1})\right) . 
\end{equation*}%
Hence, we obtain a sequence $\{x_{n}\}$ in $X$ such that for $x_{n}\in
T_{i}(x_{n-1})$, there exist $x_{n+1}\in T_{i+1}\left( x_{n}\right) $ with $%
\left( x_{n},x_{n+1}\right) \in \Delta _{2}$ such that%
\begin{equation*}
\tau +F\left( d(x_{n},x_{n+1})\right) \leq F\left( d(x_{n-1},x_{n})\right) . 
\end{equation*}%
Therefore%
\begin{eqnarray}
F\left( d(x_{n},x_{n+1})\right) &\leq &F\left( d(x_{n-1},x_{n})\right) -\tau
\leq F\left( d(x_{n-2},x_{n-1})\right) -2\tau  \notag \\
&\leq &...\leq F\left( d(x_{0},x_{1})\right) -n\tau .  \TCItag{2.2}
\end{eqnarray}%
From (2.2), we obtain $\lim\limits_{n\rightarrow \infty }F\left(
d(x_{n},x_{n+1})\right) =-\infty $ that together with ($F_{2}$) gives%
\begin{equation*}
\lim\limits_{n\rightarrow \infty }d(x_{n},x_{n+1})=0. 
\end{equation*}%
From ($F_{3}$), there exists $\lambda \in \left( 0,1\right) $ such that%
\begin{equation*}
\lim\limits_{n\rightarrow \infty }[d(x_{n},x_{n+1})]^{\lambda }F\left(
d(x_{n},x_{n+1})\right) =0. 
\end{equation*}%
From (2.2) we have%
\begin{eqnarray*}
&&[d(x_{n},x_{n+1})]^{\lambda }F\left( d(x_{n},x_{n+1})\right)
-[d(x_{n},x_{n+1})]^{\lambda }F\left( d(x_{0},x_{n+1})\right) \\
&\leq &-n\tau \lbrack d(x_{n},x_{n+1})]^{\lambda }\leq 0.
\end{eqnarray*}%
On taking limit as $n\rightarrow \infty $ we obtain%
\begin{equation*}
\lim\limits_{n\rightarrow \infty }n[d(x_{n},x_{n+1})]^{\lambda }=0. 
\end{equation*}%
Hence $\lim\limits_{n\rightarrow \infty }n^{\frac{1}{\lambda }%
}d(x_{n},x_{n+1})=0$ and there exists $n_{1}\in 
\mathbb{N}
$ such that $n^{\frac{1}{\lambda }}d(x_{n},x_{n+1})\leq 1$ for all $n\geq
n_{1}.$ So we have%
\begin{equation*}
d(x_{n},x_{n+1})\leq \frac{1}{n^{1/\lambda }} 
\end{equation*}%
for all $n\geq n_{1}.$ Now consider $m,n\in 
\mathbb{N}
$ such that $m>n\geq n_{1}$, we have%
\begin{eqnarray*}
d\left( x_{n},x_{m}\right) &\leq &d\left( x_{n},x_{n+1}\right) +d\left(
x_{n+1},x_{n+2}\right) +...+d\left( x_{m-1},x_{m}\right) \\
&\leq &\sum_{i=n}^{\infty }\frac{1}{i^{1/\lambda }}.
\end{eqnarray*}%
By the convergence of the series $\sum_{i=1}^{\infty }\frac{1}{i^{1/\lambda }%
},$ we get $d\left( x_{n},x_{m}\right) \rightarrow 0$ as $n,m\rightarrow
\infty $. Therefore $\{x_{n}\}$ is a Cauchy sequence in $X.$ Since $X$ is
complete, there exists an element $x^{\ast }\in X$ such that $%
x_{n}\rightarrow x^{\ast }$ as $n\rightarrow \infty $.\smallskip

Now, if $T_{i}$ is upper semicontinuous for any $i\in \{1,2,...,m\}$, then
as $x_{2n}\in X,$ $x_{2n+1}\in T_{i}\left( x_{2n}\right) $ with $%
x_{2n}\rightarrow x^{\ast }$ and $x_{2n+1}\rightarrow x^{\ast }$ as $%
n\rightarrow \infty $ implies that $x^{\ast }\in T_{i}\left( x^{\ast
}\right) .$ Thus from (2), we get $x^{\ast }\in T_{1}\left( x^{\ast }\right)
=T_{2}\left( x^{\ast }\right) =...=T_{m}\left( x^{\ast }\right) $.\smallskip

\noindent Finally to prove (3),\ suppose the set $\cap _{i=1}^{m}Fix\left(
T_{i}\right) $ is a well-ordered.\ We are to show that $\cap
_{i=1}^{m}Fix\left( T_{i}\right) \ $is singleton. Assume on contrary that
there exist $u$ and $v$ such that $u,v\in \cap _{i=1}^{m}Fix\left(
T_{i}\right) $ but $u\neq v$. As $(u,v)\in \Delta _{2}$, so for $%
(u_{x},v_{y})\in \Delta _{2}$ implies%
\begin{eqnarray*}
\tau +F\left( d(u,v)\right) &\leq &F(M(u,v;u,v)) \\
&=&F(\max \{d(u,v),d(u,u),d(v,v),\dfrac{d\left( u,v\right) +d\left(
v,u\right) }{2}\}) \\
&=&F\left( d\left( u,v\right) \right) ,
\end{eqnarray*}%
a contradiction as $\tau >0$.\ Hence $u=v$. Conversely, if $\cap
_{i=1}^{m}Fix\left( T_{i}\right) $ is singleton, then it follows that $\cap
_{i=1}^{m}Fix\left( T_{i}\right) $ is a well-ordered. $\square $\newline
\newline

\noindent The following corollary extends and generalizes Theorem 4.1 of 
\cite{LatifBeg} and Theorem 3.4 of \cite{Rus}.

\noindent \textbf{Corollary 2.2.} \ \ Let $(X,d,\preceq )$ be an ordered
complete metric space and $T_{1},T_{2}:X\rightarrow P_{cl}(X)$ be two
multivalued mappings. Suppose that for every $(x,y)\in \Delta _{1}$ and $%
u_{x}\in T_{i}(x),$ there exists $u_{y}\in T_{j}(y)$ for $i,j\in \{1,2\}$
with $i\neq j$ such that, $(u_{x},u_{y})\in \Delta _{2}$ implies%
\begin{equation*}
\tau +F\left( d(u_{x},u_{y})\right) \leq F(M(x,y;u_{x},u_{y})), 
\end{equation*}%
where $\tau $\ is a positive real number and%
\begin{equation*}
M(x,y;u_{x},u_{y})=\max \{d(x,y),d(x,u_{x}),d(y,u_{y}),\dfrac{d\left(
x,u_{y}\right) +d\left( y,u_{x}\right) }{2}\}. 
\end{equation*}%
Then the following statements hold:

\begin{enumerate}
\item[(1)] $Fix(T_{i})\neq \emptyset $ for any $i\in \{1,2\}$ if and only if 
$Fix(T_{1})=Fix(T_{2})\neq \emptyset .$

\item[(2)] $Fix(T_{1})=Fix(T_{2})\neq \emptyset $ provided that $T_{1}$ or $%
T_{2}$\ is upper semicontinuous.

\item[(3)] $Fix(T_{1})\cap Fix(T_{2})$ is well-ordered if and only if $%
Fix(T_{1})\cap Fix(T_{2})$ is singleton set.
\end{enumerate}

\noindent \textbf{Example 2.3.}\ \ \ Let $X=\{x_{n}=\dfrac{n(n+1)}{2}:n\in
\{1,2,3,...\}\}$ endow with usual order $\leq .$ Let%
\begin{eqnarray*}
\Delta _{1} &=&\{(x,y):x\leq y\text{ where }x,y\in X\}\text{ and} \\
\Delta _{2} &=&\{(x,y):x<y\text{ where }x,y\in X\}.
\end{eqnarray*}%
Define $T_{1}$, $T_{2}:X\rightarrow P_{cl}(X)$ as follows:%
\begin{eqnarray*}
T_{1}\left( x\right) &=&\{x_{1}\}\text{ for }x\in X,\text{ and } \\
T_{2}\left( x\right) &=&\left\{ 
\begin{array}{lll}
\{x_{1}\} & , & \text{ }x=x_{1} \\ 
\{x_{1},x_{n-1}\} & , & \text{ }x=x_{n},\text{ for }n>1\text{.}%
\end{array}%
\right.
\end{eqnarray*}%
Take $F\left( \alpha \right) =\ln \alpha +\alpha ,$ $\alpha >0$ and $\tau
=1. $ For a Euclidean metric $d$ on $X,$ and $\left( u_{x},u_{y}\right) \in
\Delta _{2},$ we consider the following cases:

\begin{enumerate}
\item[(i)] If $x=x_{1},y=x_{m},$ for $m>1,$ then for $u_{x}=x_{1}\in
T_{1}\left( x\right) ,$ there exists $u_{y}=x_{m-1}\in T_{2}\left( y\right)
, $ such that%
\begin{eqnarray*}
d(u_{x},u_{y})e^{d(u_{x},u_{y})-M\left( x,y;u_{x},u_{y}\right) } &\leq
&d(u_{x},u_{y})e^{d(u_{x},u_{y})-d\left( x,y\right) } \\
&=&\frac{m^{2}-m-2}{2}e^{-m} \\
&<&\frac{m^{2}+m-2}{2}e^{-1} \\
&=&e^{-1}d\left( x,y\right) \\
&\leq &e^{-1}M\left( x,y;u_{x},u_{y}\right) .
\end{eqnarray*}

\item[(ii)] If $x=x_{n},$ $y=x_{n+1}$ with $n>1,$ then for $u_{x}=x_{1}\in
T_{1}\left( x\right) ,$ there exists $u_{y}=x_{n-1}\in T_{2}\left( y\right)
, $ such that%
\begin{eqnarray*}
d(u_{x},u_{y})e^{d(u_{x},u_{y})-M\left( x,y;u_{x},u_{y}\right) } &\leq
&d(u_{x},u_{y})e^{d(u_{x},u_{y})-[\tfrac{d\left( x,u_{y}\right) +d\left(
y,u_{x}\right) }{2}]} \\
&=&\frac{n^{2}-n-2}{2}e^{\frac{-3n-2}{2}} \\
&<&\frac{n^{2}+4n}{2}e^{-1} \\
&=&e^{-1}[\dfrac{d\left( x,u_{y}\right) +d\left( y,u_{x}\right) }{2}] \\
&\leq &e^{-1}M_{1}\left( x,y;u_{x},u_{y}\right) .
\end{eqnarray*}

\item[(iii)] When $x=x_{n},$ $y=x_{m}$ with $m>n>1,$ then for $%
u_{x}=x_{1}\in T_{1}\left( x\right) ,$ there exists $u_{y}=x_{n-1}\in
T_{2}\left( y\right) ,$ such that%
\begin{eqnarray*}
d(u_{x},u_{y})e^{d(u_{x},u_{y})-M\left( x,y;u_{x},u_{y}\right) } &\leq
&d(u_{x},u_{y})e^{d(u_{x},u_{y})-d\left( x,u_{x}\right) } \\
&=&\frac{n^{2}-n-2}{2}e^{-n} \\
&<&\frac{n^{2}+n-2}{2}e^{-1} \\
&=&e^{-1}d\left( x,u_{x}\right) \\
&\leq &e^{-1}M\left( x,y;u_{x},u_{y}\right) .
\end{eqnarray*}
\end{enumerate}

\noindent Now we show that for $x,y\in X$, $u_{x}\in T_{2}\left( x\right) $;
there exists $u_{y}\in T_{1}\left( y\right) $ such that $\left(
u_{x},u_{y}\right) \in \Delta _{2}$ and (2.1) is satisfied. For this, we
consider the following cases:

\begin{enumerate}
\item[(i)] If $x=x_{n},$ $y=x_{1}$ with $n>1,$ we have for $u_{x}=x_{n-1}\in
T_{2}\left( x\right) ,$ there exists $u_{y}=x_{1}\in T_{1}\left( y\right) ,$
such that%
\begin{eqnarray*}
d(u_{x},u_{y})e^{d(u_{x},u_{y})-M\left( x,y;u_{x},u_{y}\right) } &\leq
&d(u_{x},u_{y})e^{d(u_{x},u_{y})-d\left( x,y\right) } \\
&=&\frac{n^{2}-n-2}{2}e^{-n} \\
&<&\frac{n^{2}+n-2}{2}e^{-1} \\
&=&e^{-1}d\left( x,y\right) \\
&\leq &e^{-1}M\left( x,y;u_{x},u_{y}\right) .
\end{eqnarray*}

\item[(ii)] In case $x=x_{n},$ $y=x_{m}$ with $m>n>1,$ then for $%
u_{x}=x_{n-1}\in T_{2}\left( x\right) ,$ there exists $u_{y}=x_{1}\in
T_{2}\left( y\right) ,$ such that%
\begin{eqnarray*}
d(u_{x},u_{y})e^{d(u_{x},u_{y})-M\left( x,y;u_{x},u_{y}\right) } &\leq
&d(u_{x},u_{y})e^{d(u_{x},u_{y})-d\left( y,u_{y}\right) } \\
&=&\frac{n^{2}-n-2}{2}e^{n^{2}-n-m^{2}-m} \\
&<&\frac{m^{2}+m-2}{2}e^{-1} \\
&=&e^{-1}d\left( y,u_{y}\right) \\
&\leq &e^{-1}M\left( x,y;u_{x},u_{y}\right) .
\end{eqnarray*}
\end{enumerate}

\noindent Thus for all $x,y$ in $X$, (2.1) is satisfied. Hence all the
conditions of Theorem 2.1 are satisfied. Moreover, $x_{1}=1$ is the unique
common fixed point of $T_{1}$ and $T_{2}$ with $Fix(T_{1})=Fix(T_{2}).$ $%
\square $\newline

\noindent The following results generalizes Theorem 3.4 of \cite{Rus} and
Theorem 3.4 of \cite{Sgr}.

\noindent \textbf{Theorem 2.4.} \ \ Let $(X,d,\preceq )$ be an ordered
complete metric space and $\{T_{i}\}_{i=1}^{m}:X\rightarrow P_{cl}(X)$ be
family of multivalued mappings. Suppose that for every $(x,y)\in \Delta _{1}$
and $u_{x}\in T_{i}(x),$ there exists $u_{y}\in T_{i+1}(y)$ for $i\in
\{1,2,...,m\}$ (with $T_{m+1}=T_{1}$ by convention) such that, $%
(u_{x},u_{y})\in \Delta _{2}$ implies%
\begin{equation}
\tau +F\left( d(u_{x},u_{y})\right) \leq F(M_{2}(x,y;u_{x},u_{y})), 
\tag{2.3}
\end{equation}%
where $\tau $\ is a positive real number and%
\begin{equation*}
M_{2}(x,y;u_{x},u_{y})=\alpha d(x,y)+\beta d(x,u_{x})+\gamma
d(y,u_{y})+\delta _{1}d\left( x,u_{y}\right) +\delta _{2}d\left(
y,u_{x}\right) , 
\end{equation*}%
and $\alpha ,\beta ,\gamma ,\delta _{1},\delta _{2}\geq 0,$ $\delta _{1}\leq
\delta _{2}$ with $\alpha +\beta +\gamma +\delta _{1}+\delta _{2}\leq 1$.
Then the following statements hold:

\begin{enumerate}
\item[(1)] $Fix(T_{i})\neq \emptyset $ for any $i\in \{1,2,...,m\}$ if and
only if $Fix(T_{1})=Fix(T_{2})=...=Fix(T_{m})\neq \emptyset .$

\item[(2)] $Fix(T_{1})=Fix(T_{2})=...=Fix(T_{m})\neq \emptyset $ provided
that any one $T_{i}$ for $i\in \{1,2,...,m\}$\ is upper semicontinuous.

\item[(3)] $\cap _{i=1}^{m}Fix(T_{i})$ is well-ordered if and only if $\cap
_{i=1}^{m}Fix(T_{i})$ is singleton set.
\end{enumerate}

\noindent \textbf{Proof.} \ \ To prove (1), let $x^{\ast }\in T_{k}(x^{\ast
})$ for some $k\in \{1,2,...,m\}.$ Assume that $x^{\ast }\notin
T_{k+1}\left( x^{\ast }\right) ,$ then there exists an $x\in T_{k+1}\left(
x^{\ast }\right) $ with $\left( x^{\ast },x\right) \in \Delta _{2}\ $such
that%
\begin{equation*}
\tau +F\left( d(x^{\ast },x)\right) \leq F(M_{2}(x^{\ast },x^{\ast };x^{\ast
},x)), 
\end{equation*}%
where%
\begin{eqnarray*}
M_{2}(x^{\ast },x^{\ast };x^{\ast },x) &=&\alpha d(x^{\ast },x^{\ast
})+\beta d(x^{\ast },x^{\ast })+\gamma d(x,x^{\ast }) \\
&&+\delta _{1}d(x^{\ast },x)+\delta _{2}d(x^{\ast },x^{\ast }) \\
&=&(\gamma +\delta _{1})d(x,x^{\ast }),
\end{eqnarray*}%
implies that%
\begin{eqnarray*}
\tau +F\left( d(x^{\ast },x)\right) &\leq &F((\gamma +\delta _{1})d(x^{\ast
},x)) \\
&\leq &F(d(x^{\ast },x)),
\end{eqnarray*}%
a contradiction as $\tau >0$. Thus $x^{\ast }=x$. Thus $x^{\ast }\in
T_{k+1}\left( x^{\ast }\right) $ and so $Fix(T_{k})\subseteq Fix(T_{k+1}).$
Similarly, we obtain that\ $Fix(T_{k+1})\subseteq Fix(T_{k+2})$ and
continuing this way, we get $Fix(T_{1})=Fix(T_{2})=...Fix(T_{k}).$ The
converse is straightforward.\smallskip

\noindent To prove (2), suppose that $x_{0}$ is an arbitrary point of $X.$
If $x_{0}\in T_{k_{0}}\left( x_{0}\right) $ for any $k_{0}\in \{1,2,...,m\},$
then by using (1), the proof is finishes. So we assume that $x_{0}\notin
T_{k_{0}}\left( x_{0}\right) $ for any $k_{0}\in \{1,2,...,m\}.$ Now for $%
i\in \{1,2,...,m\}$, if $x_{1}\in T_{i}(x_{0}),$ then there exists $x_{2}\in
T_{i+1}(x_{1})$ with $(x_{1},x_{2})\in \Delta _{2}$ such that%
\begin{equation*}
\tau +F\left( d(x_{1},x_{2})\right) \leq F(M_{2}(x_{0},x_{1};x_{1},x_{2})), 
\end{equation*}%
where%
\begin{eqnarray*}
M_{2}(x_{0},x_{1};x_{1},x_{2}) &=&\alpha d(x_{0},x_{1})+\beta
d(x_{0},x_{1})+\gamma d(x_{1},x_{2}) \\
&&+\delta _{1}d(x_{0},x_{2})+\delta _{2}d(x_{1},x_{1}) \\
&\leq &(\alpha +\beta +\delta _{1})d(x_{0},x_{1})+(\gamma +\delta
_{1})d(x_{1},x_{2}).
\end{eqnarray*}%
Now, if $d(x_{0},x_{1})\leq d(x_{1},x_{2}),$ then we have%
\begin{eqnarray*}
\tau +F\left( d(x_{1},x_{2})\right) &\leq &F((\alpha +\beta +\gamma +2\delta
_{1})d(x_{1},x_{2})) \\
&\leq &F(d(x_{1},x_{2})),
\end{eqnarray*}%
a contradiction. Therefore%
\begin{equation*}
\tau +F\left( d(x_{1},x_{2})\right) \leq F\left( d(x_{0},x_{1})\right) . 
\end{equation*}%
Next for this $x_{2}\in T_{i+1}\left( x_{1}\right) ,$ there exists $x_{3}\in
T_{i+2}(x_{2})$ with $\left( x_{2},x_{3}\right) \in \Delta _{2}$ such that%
\begin{equation*}
\tau +F\left( d(x_{2},x_{3})\right) \leq F(M_{2}(x_{1},x_{2};x_{2},x_{3})), 
\end{equation*}%
where%
\begin{eqnarray*}
M_{2}(x_{1},x_{2};x_{2},x_{3}) &=&\alpha d(x_{1},x_{2})+\beta
d(x_{1},x_{2})+\gamma d(x_{2},x_{3}) \\
&&+\delta _{1}d(x_{1},x_{3})+\delta _{2}d(x_{2},x_{2}) \\
&\leq &(\alpha +\beta +\delta _{1})d(x_{1},x_{2})+(\gamma +\delta
_{1})d(x_{2},x_{3}).
\end{eqnarray*}%
Now, if $d(x_{1},x_{2})\leq d(x_{2},x_{3})$ then%
\begin{eqnarray*}
\tau +F\left( d(x_{2},x_{3})\right) &\leq &F((\alpha +\beta +\gamma +2\delta
_{1})d(x_{2},x_{3})) \\
&\leq &F\left( d\left( x_{2},x_{3}\right) \right) ,
\end{eqnarray*}%
a contradiction as $\tau >0$. Therefore%
\begin{equation*}
\tau +F\left( d(x_{2},x_{3})\right) \leq F\left( d(x_{1},x_{2})\right) . 
\end{equation*}%
Continuing this process, for $x_{2n}\in T_{i}(x_{2n-1})$, there exist $%
x_{2n+1}\in T_{i+1}\left( x_{2n}\right) $ with $\left(
x_{2n},x_{2n+1}\right) \in \Delta _{2}$ such that%
\begin{equation*}
\tau +F\left( d(x_{2n},x_{2n+1})\right) \leq F\left(
M_{2}(x_{2n-1},x_{2n};x_{2n},x_{2n+1})\right) , 
\end{equation*}%
where%
\begin{eqnarray*}
M_{2}(x_{2n-1},x_{2n};x_{2n},x_{2n+1}) &=&\alpha d(x_{2n-1},x_{2n})+\beta
d(x_{2n-1},x_{2n})+\gamma d(x_{2n},x_{2n+1}) \\
&&+\delta _{1}d(x_{2n-1},x_{2n+1})+\delta _{2}d(x_{2n},x_{2n}) \\
&\leq &\left( \alpha +\beta +\delta _{1}\right) d(x_{2n-1},x_{2n})+\left(
\gamma +\delta _{1}\right) d(x_{2n},x_{2n+1}) \\
&\leq &d(x_{2n-1},x_{2n}),
\end{eqnarray*}%
that is,%
\begin{equation*}
\tau +F\left( d(x_{2n},x_{2n+1})\right) \leq F\left(
d(x_{2n-1},x_{2n})\right) . 
\end{equation*}%
Similarly, for $x_{2n+1}\in T_{i+1}(x_{2n})$, there exist $x_{2n+2}\in
T_{i+2}\left( x_{2n+1}\right) $ such that for $\left(
x_{2n+1},x_{2n+2}\right) \in \Delta _{2}$ implies%
\begin{equation*}
\tau +F\left( d(x_{2n+1},x_{2n+2})\right) \leq F\left(
d(x_{2n},x_{2n+1})\right) . 
\end{equation*}%
Hence, we obtain a sequence $\{x_{n}\}$ in $X$ such that for $x_{n}\in
T_{i}(x_{n-1})$, there exist $x_{n+1}\in T_{i+1}\left( x_{n}\right) $ with $%
\left( x_{n},x_{n+1}\right) \in \Delta _{2}$ such that%
\begin{equation*}
\tau +F\left( d(x_{n},x_{n+1})\right) \leq F\left( d(x_{n-1},x_{n})\right) . 
\end{equation*}%
Therefore%
\begin{eqnarray}
F\left( d(x_{n},x_{n+1})\right) &\leq &F\left( d(x_{n-1},x_{n})\right) -\tau
\leq F\left( d(x_{n-2},x_{n-1})\right) -2\tau  \notag \\
&\leq &...\leq F\left( d(x_{0},x_{1})\right) -n\tau .  \TCItag{2.4}
\end{eqnarray}%
From (2.4), we obtain $\lim\limits_{n\rightarrow \infty }F\left(
d(x_{n},x_{n+1})\right) =-\infty $ that together with ($F_{2}$) gives%
\begin{equation*}
\lim\limits_{n\rightarrow \infty }d(x_{n},x_{n+1})=0. 
\end{equation*}%
Follows the arguments those in proof of Theorem 2.1, $\{x_{n}\}$ is a Cauchy
sequence in $X.$ Since $X$ is complete, there exists an element $x^{\ast
}\in X$ such that $x_{n}\rightarrow x^{\ast }$ as $n\rightarrow \infty $%
.\smallskip

Now, if $T_{i}$ is upper semicontinuous for any $i\in \{1,2,...,m\}$, then
as $x_{2n}\in X,$ $x_{2n+1}\in T_{i}\left( x_{2n}\right) $ with $%
x_{2n}\rightarrow x^{\ast }$ and $x_{2n+1}\rightarrow x^{\ast }$ as $%
n\rightarrow \infty $ implies that $x^{\ast }\in T_{i}\left( x^{\ast
}\right) .$ Thus from (2), we get $x^{\ast }\in T_{1}\left( x^{\ast }\right)
=T_{2}\left( x^{\ast }\right) =...=T_{m}\left( x^{\ast }\right) $.\smallskip

\noindent Finally to prove (3),\ suppose the set $\cap _{i=1}^{m}Fix\left(
T_{i}\right) $ is a well-ordered.\ We are to show that $\cap
_{i=1}^{m}Fix\left( T_{i}\right) \ $is singleton. Assume on contrary that
there exist $u$ and $v$ such that $u,v\in \cap _{i=1}^{m}Fix\left(
T_{i}\right) $ but $u\neq v$. As $(u,v)\in \Delta _{2}$, so for $%
(u_{x},v_{y})\in \Delta _{2}$ implies%
\begin{equation*}
\tau +F\left( d(u,v)\right) \leq F(M_{2}(u,v;u,v)), 
\end{equation*}%
where%
\begin{eqnarray*}
M_{2}(u,v;u,v) &=&\alpha d(u,v)+\beta d(u,u)+\gamma d(v,v) \\
&&+\delta _{1}d\left( u,v\right) +\delta _{2}d\left( v,u\right) \\
&=&\left( \alpha +\delta _{1}+\delta _{2}\right) d\left( x,y\right) ,
\end{eqnarray*}%
that is,%
\begin{eqnarray*}
\tau +F\left( d(u,v)\right) &\leq &F\left( \left( \alpha +\delta _{1}+\delta
_{2}\right) d\left( x,y\right) \right) \\
&\leq &F\left( d\left( u,v\right) \right) ,
\end{eqnarray*}%
a contradiction as $\tau >0$.\ Hence $u=v$. Conversely, if $\cap
_{i=1}^{m}Fix\left( T_{i}\right) $ is singleton, then it follows that $\cap
_{i=1}^{m}Fix\left( T_{i}\right) $ is a well-ordered. $\square $\newline

\noindent The following corollary extends Theorem 3.1 of \cite{Rus}, in the
case of family of mappings in ordered metric space.

\noindent \textbf{Corollary 2.5.} \ \ Let $(X,d,\preceq )$ be an ordered
complete metric space and $\{T_{i}\}_{i=1}^{m}:X\rightarrow P_{cl}(X)$ be
family of multivalued mappings. Suppose that for every $(x,y)\in \Delta _{1}$
and $u_{x}\in T_{i}(x),$ there exists $u_{y}\in T_{i+1}(y)$ for $i\in
\{1,2,...,m\}$ (with $T_{m+1}=T_{1}$ by convention) such that, $%
(u_{x},u_{y})\in \Delta _{2}$ implies%
\begin{equation}
\tau +F\left( d(u_{x},u_{y})\right) \leq F(\alpha d\left( x,y\right) +\beta
d(x,u_{x})+\gamma d(y,u_{y})]),  \tag{2.5}
\end{equation}%
where $\tau $\ is a positive real number and $\alpha ,\beta ,\gamma \geq 0$
with $\alpha ,\beta ,\gamma \leq 1.$ Then the conclusions obtained in
Theorem 2.3 remain true.\bigskip

\noindent The following corollary extends Theorem 4.1 of \cite{LatifBeg}.

\noindent \textbf{Corollary 2.6.} \ \ Let $(X,d,\preceq )$ be an ordered
complete metric space and $\{T_{i}\}_{i=1}^{m}:X\rightarrow P_{cl}(X)$ be
family of multivalued mappings. Suppose that for every $(x,y)\in \Delta _{1}$
and $u_{x}\in T_{i}(x),$ there exists $u_{y}\in T_{i+1}(y)$ for $i\in
\{1,2,...,m\}$ (with $T_{m+1}=T_{1}$ by convention) such that, $%
(u_{x},u_{y})\in \Delta _{2}$ implies%
\begin{equation}
\tau +F\left( d(u_{x},u_{y})\right) \leq F(h[d(x,u_{x})+d(y,u_{y})]), 
\tag{2.6}
\end{equation}%
where $\tau $\ is a positive real number and $h\in \lbrack 0,\dfrac{1}{2}].$
Then the conclusions obtained in Theorem 2.3 remain true.\bigskip

\noindent \textbf{Corollary 2.7.} \ \ Let $(X,d,\preceq )$ be an ordered
complete metric space and $\{T_{i}\}_{i=1}^{m}:X\rightarrow P_{cl}(X)$ be
family of multivalued mappings. Suppose that for every $(x,y)\in \Delta _{1}$
and $u_{x}\in T_{i}(x),$ there exists $u_{y}\in T_{i+1}(y)$ for $i\in
\{1,2,...,m\}$ (with $T_{m+1}=T_{1}$ by convention) such that, $%
(u_{x},u_{y})\in \Delta _{2}$ implies%
\begin{equation}
\tau +F\left( d(u_{x},u_{y})\right) \leq F(d(x,y)),  \tag{2.7}
\end{equation}%
where $\tau $\ is a positive real number. Then the conclusions obtained in
Theorem 2.3 remain true.\bigskip

\noindent The above corollary extends Theorem 4.1 of \cite{LatifBeg}.\bigskip

\noindent \textbf{Conclusion.} Recently many results appeared in the
literature giving the problems related to the common fixed point for
multivalued maps. I this paper we obtained the results for existence of
common fixed points of family of maps that satisfying generalized $F$%
-contractions in ordered structured metric spaces. We presented some
examples to show the validity of established results.

\noindent \textbf{Acknowledgement.} The first author is grateful to the
Erasmus Mundus project FUSION for the postdoctoral fellowship visiting to M%
\"{a}lardalen University Sweden and to the Division of Applied Mathematics
at the School of Education, Culture and Communication for creating excellent
research environment.

\end{document}